\documentclass[a4paper, draft, reqno]{compositio}
\usepackage{epsf}
\usepackage{amssymb}
\usepackage{amsmath}
\begin{document}
%{{{ macros

\title[Towers of recollement]
{Representation theory of 
towers of recollement:\\ theory, notes, and examples}

\author{Anton Cox}
\email{A.G.Cox@city.ac.uk}
\address{Centre for Mathematical Science\\
  City University\\ Northampton Square\\ London\\ EC1V
  0HB\\ UK}  

\author{Paul Martin}
\email{P.P.Martin@city.ac.uk}
\address{Centre for Mathematical Science\\
  City University\\ Northampton Square\\ London\\ EC1V
  0HB\\ UK}  

\author{Alison Parker}
\email{alisonp@maths.usyd.edu.au}
\address{School of Mathematics and Statistics F07\\
University of Sydney\\ NSW 2006\\ Australia}

\author{Changchang Xi}
\email{xicc@bnu.edu.cn}
\address{Department of Mathematics\\ Beijing Normal
University\\  100875 Beijing\\ China}

\shortauthors{A. Cox et al.}
\classification{20C08}
\thanks{The first author is supported by Nuffield grant scheme
NUF-NAL 02,  and the fourth by a China-UK
joint project of the Royal Society (No. 15262)}

\begin{abstract}
We give an axiomatic framework for studying the representation theory
of towers of algebras. We introduce a new class of algebras, contour
algebras, generalising (and interpolating between) blob algebras and
cyclotomic Temperley-Lieb algebras. We demonstrate the utility of our
formalism by applying it to this class.
\end{abstract}
\maketitle

\theoremstyle{plain}
\newtheorem{Theorem}{Theorem}
\newtheorem{thm}{Theorem}[section]
\newtheorem*{introthm}{Theorem}
\newtheorem{prop}[thm]{Proposition}
\newtheorem{lem}[thm]{Lemma}
\newtheorem{cor}[thm]{Corollary}
\newtheorem{con}[thm]{Conjecture}
\theoremstyle{remark}
\newtheorem{rem}[thm]{Remark}
\newtheorem{ex}[thm]{Examples}

\newcommand{\sgn}[0]{\mbox{\rm sgn }}
\newcommand{\zed}{{\mathbb Z}}
\newcommand{\NN}{{\mathbb N}}
\newcommand{\CC}{{\mathbb C}}
\newcommand{\pf}[0]{{\noindent\bf Proof:} }
\newcommand{\qdet}[0]{q\mbox{\rm -det}}
\newcommand{\Hom}[0]{\mbox{\rm Hom}}
\newcommand{\Ext}[0]{\mbox{\rm Ext}}
\newcommand{\End}[0]{\mbox{\rm End}}
\newcommand{\soc}[0]{\mbox{\rm soc}}
\newcommand{\hd}[0]{\mbox{\rm hd}}
\newcommand{\ch}[0]{\mbox{\rm ch}}
\newcommand{\dual}[0]{^{\ast}}
\newcommand{\frob}[0]{^{\mbox{\rm\tiny F}}}
\newcommand{\tlam}[0]{\tilde{\lambda}}
\newcommand{\bmu}[0]{\bar{\mu}}
\newcommand{\too}[0]{\rightarrow}
\newcommand{\tooo}[0]{\longrightarrow}
\newcommand{\Exti}[0]{\Ext^1_{G_1}}
\newcommand{\Homi}[0]{\Hom_{G_1}}
\newcommand{\otherwise}[0]{\mbox{\rm otherwise}}
\newcommand{\wif}[0]{\mbox{\rm if }}
\newcommand{\wand}[0]{\mbox{\rm and }}
\renewcommand{\mod}[0]{\mbox{\rm mod }}
\newcommand{\congG}[0]{\cong_{_{G_1}}}
\newcommand{\even}[0]{\ \mbox{\rm even}}
\newcommand{\odd}[0]{\ \mbox{\rm odd}}
\newcommand{\Extg}[0]{\Ext^1_G}
\newcommand{\bg}[0]{\bar{G}}
\newcommand{\gi}[0]{G_1}
\newcommand{\gl}[0]{\mbox{\rm GL}_2}
\newcommand{\SL}[0]{\mbox{\rm SL}_2}
\newcommand{\tgl}[0]{\mbox{\tiny\rm GL}_2}
\newcommand{\tsl}[0]{\mbox{\tiny\rm SL}_2}
\newcommand{\Homgl}[0]{\Hom_{\mbox{\tiny\rm GL}_2}}
\newcommand{\upgl}[0]{^{\mbox{\tiny\rm GL}_2}}
\newcommand{\Extgl}[0]{\Ext^1_{\tgl}}
\newcommand{\Extsl}[0]{\Ext^1_{\tsl}}
\newcommand{\ha}[0]{\hat{a}}
\newcommand{\hr}[0]{\hat{r}}
\newcommand{\wfor}[0]{\mbox{\rm for }}
\newcommand{\mrd}[0]{M_rD}
\newcommand{\ind}[0]{\mbox{\rm ind}}
\newcommand{\res}[0]{\mbox{\rm res}}
\newcommand{\supp}[0]{\mbox{\rm supp}}
\newcommand{\wt}[0]{\mbox{\rm wt}}
\newcommand{\fit}[0]{\phantom{\Big(}}

%}}}

%{{{ Authors

%}}}
%{{{ Intro
\section*{Introduction}

Let $A $ be a finite-dimensional algebra and $e \in A $ be an
idempotent.  The category $eAe $-mod is fully embedded in $A $-mod and
the remaining simples $L $ for $A $ are characterised by $eL =0 $.  In
particular, we have an exact `localisation' functor
$$\begin{array}{ccc} F:A \mbox{\rm -mod} &\too& eAe \mbox{\rm -mod}\\
M & \longmapsto &eM \end{array} $$ 
which takes simples to simples or
zero.  Indeed, every simple $eAe $-module arises in this way:
\medskip\smallskip

\begin{introthm}[(Green \cite{green})]
Let $\{L (\lambda): \lambda \in
\Lambda \} $ be a full set 
of simple $A $-modules, and set $\Lambda^e = \{ \lambda \in \Lambda:eL(
\lambda) \neq0 \} $.  Then $\{eL( \lambda): \lambda \in \Lambda^e \} $
is a full set of simple $eAe $-modules.
Further, the simple modules $L (\lambda) $ with $\lambda \in
\Lambda\backslash \Lambda^e $ are a full set of simple $A/AeA
$-modules.
\end{introthm}

We define the globalisation functor by 
$$G :N \longmapsto Ae \otimes_{eAe }N $$ and note that $FG (N) \cong
N$ and $G$ is a full embedding.  Cline, Parshall and Scott \cite{cps1}
use this idea to provide examples of recollement \cite{bbd} in the
context of quasi-heredity and highest weight categories.  Following an
application to the Temperley-Lieb algebra in \cite{marbook}, the
second author and Saleur then used it for the tower $b_1 \subset b_2
\subset \ldots $ of blob algebras \cite{msblob}, for which there exist
idempotents $e_n \in b_n $ such that $e_nb_ne_n \cong b_{n-2} $, to
recursively analyse the representation theory of the entire tower.

There are in fact a significant number of interesting towers of
algebras with such idempotents, particularly among algebras equipped
with a diagram calculus and algebras arising in invariant theory
\cite{brauer37, brown55, marbook, mardef,jonesplan, greentangle,
mazbrauer, mazend, bloss, marel, rx}.  In Section \ref{towers} we
abstract and formalise aspects of the common procedure used to analyse
such towers of algebras in \cite{msblob,mardef} (while largely
avoiding the explicit construction of bases).  In Sections \ref{rxsec}
and \ref{rxrep} we demonstrate the utility of this formalism by
applying it to a new class of diagram algebras, the contour
algebras. This is a collection of towers of algebras which includes as
special cases the Temperley-Lieb algebras and blob algebras, and the
cyclotomic Temperley-Lieb algebras recently defined by Rui and the
last author \cite{rx}. The formalism allows us to index simple modules
very easily, to construct standard modules, and to locate many
standard module morphisms efficiently.  In Section \ref{gramres} we
carry out the algebra-specific calculations required by our
formalism. Finally in Section \ref{discus} we return to a discussion
of our axiom scheme. We explore the consequences of modifying our
axioms at various points, and the relationship between them and other
such exercises in the literature.

Our notion of a tower of recollement combines certain ideas from the tower
formalism in \cite{ghj} (but relaxing the emphasis on semisimplicity)
with the notion of recollement in \cite{cps1}. (The latter is a special
case of the general notion of recollement in \cite{bbd}.) We only make
explicit one of the two defining functors in a recollement diagram; the
other is implicit in this approach (see \cite[Section 2]{cps1}) but
not needed in this paper.

Although we will make no use of it in what follows, it is worth
remarking on the physics that originally drove this approach. These
algebras (over $\CC$) are transfer matrix algebras in the sense of
\cite{marbook}. The physical context naturally brings two properties
into play. First that the algebras arise as a tower (corresponding to
different physical system sizes), and second that their module
categories embed in each other (corresponding to the thermodynamic
limit). It is the interplay between these two ways of passing through
the tower that lies at the heart of our axiomatisation.

%}}}
%{{{ S Towers of recollement
\section{Towers of recollement}\label{towers}
%{{{ axioms

Let $A_n $ (with $n \geq 0 $) be a family of finite-dimensional
algebras, with idempotents $e_n$ in $A_n $.  For simplicity we shall
assume that $A_n$ is defined over an algebraically closed field
$k$. We will impose a series of restrictions on such algebras
sufficient for an analysis of their representation theory along the
lines of that carried out in \cite{msblob}.  The rationale for
introducing axioms (A1--6), which now follow, is that they allow us to
inductively classify the simple $A_n$-modules, and to determine which
of the algebras in the family are semisimple (along with lots of
homological data when they are not), with only a minimum of
calculations.

We first assume

\medskip
\newcounter{axiom}
\setcounter{axiom}{1}
\begin{list}{\quad \quad{\bf (A\arabic{axiom})}}
\item  {\it For each $n \geq 2 $ we have an isomorphism 
$$ \Phi_n:A_{n-2}\longrightarrow e_nA_ne_n. $$}
\end{list}
\medskip

With this assumption we define a pair of families of functors $F_n :
A_n \mbox{\rm -mod} \too A_{n-2} \mbox{\rm -mod}$ and $G_n: A_{n}
\mbox{\rm -mod} \too A_{n+2} \mbox{\rm -mod}$ as  in the
introduction.  That is, $F_n(M)= e_nM$ and
$G_{n-2}(N)=A_{n}e_{n}\otimes_{e_nA_ne_n}N$ (where in each case we are
using the isomorphism in (A1)).  Note that the right inverse to $F_n$
is $G_{n-2}$.  

Denote the indexing set for the simple $A_n$-modules by
$\Lambda_n$, and that for the simple $A_n/A_ne_nA_n$-modules by
$\Lambda^n$.  Then by (A1) and the Theorem in the introduction we have
\begin{equation}\label{simplesare}
\Lambda_n=\Lambda^n\sqcup\Lambda_{n-2}
\end{equation}
and hence,  provided that $\Lambda_0$, $\Lambda_1$ and $\Lambda^n$ are known,
%$A_0$, $A_1$, and the successive
%quotients $A_n/A_ne_nA_n$ have a known structure, 
this immediately
allows the simple modules for each $A_n$ to be classified by
induction. We will illustrate this by providing a very short proof of
the classification of simple modules for the contour algebras in
Corollary \ref{classsimple}.

By (\ref{simplesare}) we may regard $\Lambda_n$ as a subset of
$\Lambda_{n+2}$, and set
$\Lambda=(\lim_n\Lambda_{2n})\sqcup(\lim_n\Lambda_{2n+1})$. We call
elements of $\Lambda$ {\it weights}.  For $m,n\in\NN$ with $m-n$ even
we set $\Lambda^n_m=\Lambda^n$ {\it regarded as a subset of
$\Lambda_m$} if $m\geq n$, and $\Lambda^n_m=\emptyset$ otherwise.
%For each $j$ we
%identify $A_{j-2}$ with $e_jA_je_j$ via the isomorphism in (A1). 

Set $e_{n,0}=1$ in $A_n$, and for $1\leq i\leq \frac{n}{2}$ 
define new idempotents in $A_n$ by setting $e_{n,i}=\Phi_n(e_{n-2,i-1})$.
 To these elements we associate corresponding
quotients of $A_n$ by setting $A_{n,i}=A_n/(A_ne_{n,i+1}A_n)$.

It will be convenient to
have the machinery of quasi-heredity at our disposal. For this reason
we next assume

\medskip
\begin{list}{\quad \quad{\bf (A\arabic{axiom})}}{\usecounter{axiom}}
\setcounter{axiom}{1}
\item {\it (i) The algebra $A_n/A_ne_nA_n$ is semisimple.

(ii) For each $n\geq 0$ and $0\leq i\leq \frac{n}{2}$, setting
$e=e_{n,i}$ and $A=A_{n,i}$, the surjective multiplication map
$Ae\otimes_{eAe}eA\too
  AeA$ is a bijection.}
%\item For each $n\geq 0$ and $0\leq i\leq \frac{n}{2}$ we have
%
%(i) $e_{n,i} A_{n,i} e_{n,i}$ is semisimple.
%
%(ii) the surjective multiplication map
%$A_{n,i}e_{n,i}\otimes_{e_{n,i}A_{n,i}e_{n,i}}e_{n,i}A_{n,i}\too
%  A_{n,i}e_{n,i}A_{n,i}$ is a bijection.
\end{list}
\medskip

Note that condition (i) (with (A1)) implies that $e_{n,i}A_{n,i}e_{n,i}$ is
semisimple for all $n\geq 0$ and $0\leq i\leq \frac{n}{2}$.  We have
chosen to state (A2) in the form above to emphasise the elementary
nature of the condition (and because this is the form in which it will
be verified, which is an entirely routine matter in
specific algebras, as we will exemplify in Proposition~\ref{qhprop}).
However, by \cite[Statement 7]{dr2} (or \cite[Definition 3.3.1 and the
remarks following]{martin}), it is straightforward to verify that we
could replace (A2) by

\medskip
\begin{list}{\quad \quad{\bf (A\arabic{axiom}$'$)}}{\usecounter{axiom}}
\setcounter{axiom}{1}
\item {\it  For each $n\geq 0$ the algebra $A_n$ is quasi-hereditary, with 
heredity chain of the form
$$
0\subset\cdots\subset A_ne_{n,i}A_n\subset\cdots\subset 
A_ne_{n,0}A_n = A_n.
$$}
\end{list}
\medskip

 As $A_n$ is
quasi-hereditary, there is for each $\lambda\in\Lambda_n$ a standard
module $\Delta_n(\lambda)$, with simple head $L_n(\lambda)$. If we set
$\lambda\leq\mu$ if either $\lambda=\mu$ or $\lambda\in\Lambda_n^r$  and
$\mu\in\Lambda_n^s$ with $r>s$, then all other composition factors of
$\Delta_n(\lambda)$ are labelled by weights $\mu$ with
$\mu<\lambda$. Note that for $\lambda\in\Lambda^n_n$, we have that
$\Delta_n(\lambda)\cong L_n(\lambda)$, and that this is just the lift
of a simple module for the quotient algebra $A_n/A_ne_nA_n$.  Arguing as in
\cite[Proposition 3]{mryom} we see that
\begin{equation}\label{Gdoes}
G_n (\Delta_n(\lambda))\cong\Delta_{n+2}(\lambda).
\end{equation}
Similarly (see for example \cite[A1]{don2}) we have
\begin{equation}
F_n(\Delta_n(\lambda))\cong\left\{\begin{array}{ll}
\Delta_{n-2}(\lambda)&\wif\lambda\in\Lambda_{n-2}\\
0&\wif \lambda\in\Lambda^n.\end{array}\right.
\end{equation}

Crucially we impose  %%require of our algebras 
a {\em second} way of passing through the family of algebras:

\medskip
\begin{list}{\quad \quad{\bf (A\arabic{axiom})}}{\usecounter{axiom}}
\setcounter{axiom}{2}
\item {\it For each $n\geq 0$ the algebra $A_n$ can be identified with a
subalgebra of $A_{n+1}$.}
\end{list}
\medskip

The other main tool we wish to use, then, is Frobenius reciprocity. 
For this  %%for which
we will need to have certain controls over induction and restriction for
our families of modules.  Essentially, we want these to have a {\em local
behaviour} and be {\em compatible} with globalisation, in a sense we
now describe. 

If a module $M$ in $ A_n$-mod has a $\Delta_n$-filtration (i.e. a
filtration with successive quotients isomorphic to some
$\Delta_n(\lambda_i)$'s) we define the {\it support} of $M$, denoted
$\supp_n(M) $, to be the set of labels $\lambda $ for which $\Delta
(\lambda) $ occurs in this filtration. (As standard modules form a
basis for the Grothendieck group of a quasi-hereditary algebra, this
is well-defined.) We will also need to consider the restriction
functor $\res_n:A_n\mbox{\rm -mod}\too A_{n-1}\mbox{\rm -mod}$ and the
induction functor $\ind_n:A_n\mbox{\rm -mod}\too A_{n+1}\mbox{\rm
-mod}$ given by $\ind_n(M)=A_{n+1}\otimes_{A_n}M$. We will omit
suffixes from $\supp_n$, $\ind_n$, $\res_n$ and $\Delta_n$-filtration
whenever this is unambiguous. Our next three assumptions ensure that
induction and restriction behave well in this setting.
\medskip

\begin{list}{\quad \quad{\bf (A\arabic{axiom})}}{\usecounter{axiom}}
\setlength{\itemindent}{0pt}
\setcounter{axiom}{3}
\item{\it  For all $n\geq 1$ we have that $A_ne_n\cong A_{n-1}$ as a left
$A_{n-1}$--, right $A_{n-2}$--bimodule.}
\end{list}
\medskip

\noindent We can immediately deduce from (A4) that
for each $\lambda \in \Lambda_n $ we have that 
\begin{equation}\label{indres}
\res (G_n (\Delta_n
(\lambda))) \cong \ind \Delta_n(\lambda).
\end{equation}

\begin{list}{\quad \quad{\bf (A\arabic{axiom})}}{\usecounter{axiom}}
\setlength{\itemindent}{0pt}
\setcounter{axiom}{4}
\item {\it For each $\lambda \in \Lambda_n^m $ we have that $\res
(\Delta_n(\lambda)) $ has a $\Delta $-filtration and 
$$\supp (\res (\Delta_n (\lambda)) \subseteq
\Lambda_{n-1}^{m-1}\sqcup
\Lambda_{n-1}^{m+1}.$$} 
\end{list}
\medskip

Equation (\ref{indres}) now implies the analogue of (A5) for
induction. Using (\ref{Gdoes}) we deduce from (A5) and (\ref{indres})
that for each $\lambda \in \Lambda_n^m $ the module $\ind
(\Delta_n(\lambda)) $ has a $\Delta $-filtration, and
\begin{equation}\label{A4b}
\supp (\ind (\Delta_n (\lambda)) \subseteq \Lambda_{n+1}^{m-1}
\sqcup
\Lambda_{n+1}^{m+1}.
\end{equation}

\begin{list}{\quad \quad{\bf (A\arabic{axiom})}}{\usecounter{axiom}}
\setlength{\itemindent}{0pt}
\setcounter{axiom}{5}
\item {\it For each $\lambda\in\Lambda_n^m$ there exists
$\mu\in\Lambda_{n-1}^{m-1}$ such that 
$\lambda\in\supp(\ind\Delta_{n-1}(\mu))$.}
\end{list}
\medskip

In the presence of (A5) this is equivalent to 

\begin{list}{\quad \quad{\bf (A\arabic{axiom}$'$)}}{\usecounter{axiom}}
\setlength{\itemindent}{0pt}
\setcounter{axiom}{5}
\item {\it For each $\lambda\in\Lambda_n^m$ there exists
$\mu\in\Lambda_{n+1}^{m-1}$ such that 
$\lambda\in\supp(\res\Delta_{n+1}(\mu))$.}
\end{list}
\medskip

For a quasi-hereditary algebra we have that
$\Ext^1(\Delta(\lambda),\Delta(\mu))\neq0$ implies that
$\lambda<\mu$. Therefore (A6) is equivalent to the requirement that
for each $\lambda\in\Lambda_n$ there exists $\mu\in\Lambda_{n-1}$ such
that there is a surjection
\begin{equation}\label{A5b}
\ind\Delta_{n-1}(\mu)\too\Delta_n(\lambda)\too 0.
\end{equation}

%One final assumption?
%
%\medskip
%\begin{list}{\quad \quad{\bf (A\arabic{axiom})}}{\usecounter{axiom}}
%\setlength{\itemindent}{0pt}
%\setcounter{axiom}{6}
%\item For each $n$ we have a simple preserving duality ${ }^o$ on mod($A$).
%\end{list}
%\medskip
%
%In fact in all the examples considered in this paper this duality
%will be induced by an isomorphism $A\cong A^{op}$, but we do not
%require this to be the case.
%
%}}}
\medskip

We shall call a family of algebras satisfying (A1--6) a {\it tower of
recollement}, since it broadly combines ideas from \cite{ghj} and
\cite{cps1} as discussed in the Introduction.

The axiomatic framework introduced so far is sufficient to reduce the
study of various general homological problems to certain explicit
calculations, as illustrated by

\begin{thm}\label{homcatch}
(i) For all pairs of weights $\lambda\in\Lambda_n^m$ and
$\mu\in\Lambda_n^l$  we have
$$\Hom(\Delta_n(\lambda),\Delta_n(\mu))\cong\left\{\begin{array}{ll}
\Hom(\Delta_m(\lambda),\Delta_m(\mu))&\wif l\leq m\\
0 &\otherwise.\end{array}\right.$$
(ii) Suppose that for all $n\geq 0$ and pairs of weights
$\lambda\in\Lambda_n^n$ and $\mu\in\Lambda_n^{n-2}$ we have 
$$\Hom(\Delta_n(\lambda),\Delta_n(\mu))=0.$$
Then each of the algebras $A_n$ is semisimple.
\end{thm}
\begin{proof} For (i) we first note that quasi-heredity implies that for any
non-zero Hom-space between standard modules as above we must have
$\lambda\leq\mu$, and hence we may assume that $l\leq m$. As each
$G_n$ is a full embedding, any non-zero homomorphism between standard
modules $\Delta_n(\lambda)$ and $\Delta_n(\mu)$ corresponds to a
morphism between some pair of standards $\Delta_m(\lambda)$ and
$\Delta_m(\mu)$ with $\lambda\in\Lambda_m^{m}$.

For (ii) we will proceed by induction on $n$. Recall that in a
quasi-hereditary algebra, the standard module $\Delta(\lambda)$ is
defined to be the largest quotient of the projective cover
$P(\lambda)$ of $L(\lambda)$ with the property that all of its
composition factors $L(\mu)$ satisfy $\mu\leq \lambda$. For
semisimplicity it is enough to show that all the $P(\lambda)$ are
simple. For any finite dimensional module $M$ we have
$\dim\Hom(P(\lambda),M)=[M:L(\lambda)]$, the multiplicity of
$L(\lambda)$ as a composition factor of $M$. Hence it is enough to
show that $\Hom(P(\lambda),P(\mu))=0$ for $\mu\neq \lambda$. As
$P(\lambda)$ has a filtration by standard modules, it is enough to
show that $\Hom(\Delta(\lambda),\Delta(\mu))=0$ for $\mu\neq\lambda$.

Suppose that $\lambda$ and $\mu$ are such that
$\Hom(\Delta_n(\lambda),\Delta_n(\mu))\neq 0$. Then in the order
induced by quasi-heredity we must have $\lambda\leq\mu$; i.e. either
$\lambda\in\Lambda_n^r$ and $\mu\in\Lambda_n^s$ with $r>s$, or
$\lambda=\mu$. In the latter case quasi-heredity implies that
$\End(\Delta_n(\lambda))\cong k$, and so we may assume that $r>s$.

If $r<n$ then $F_n\Delta_n(\lambda)\cong\Delta_{n-2}(\lambda)$ and 
$F_n\Delta_n(\mu)\cong\Delta_{n-2}(\mu)$. Further, as
$\Delta_n(\lambda)$ has simple head $L_n(\lambda)$
which is not killed by $F_n$, any non-zero
homomorphism from $\Delta_n(\lambda)$ to $\Delta_n(\mu)$ survives
under $F_n$. Hence, as $A_{n-2}$ is semisimple, there are no non-zero
morphisms between $\Delta_n(\lambda)$ and $\Delta_n(\mu)$.

Thus we may assume that $r=n$ and $s<n$. Then by (\ref{A5b}) there exists
a weight $\tau\in\Lambda_{n-1}$ such that
$\ind\Delta_{n-1}(\tau)\too\Delta_n(\lambda)\too 0$, and by (\ref{A4b}) we
have that $\tau\in\Lambda_{n-1}^{n-1}$. Now we have an injection
$$0\too\Hom(\Delta_n(\lambda),\Delta_n(\mu))\too
\Hom(\ind\Delta_{n-1}(\tau),\Delta_n(\mu))$$
and by Frobenius reciprocity we have
$$\Hom(\ind\Delta_{n-1}(\tau),\Delta_n(\mu))\cong
\Hom(\Delta_{n-1}(\tau),\res\Delta_n(\mu)).$$
By (A3) and the semisimplicity of $A_{n-2}$ we
have that 
$$\res\Delta_n(\mu)\cong
\left(\oplus_i\Delta_{n-1}(\nu_i)\right)\oplus
\left(\oplus_j\Delta(\nu_j)\right)$$ 
for some $\nu_i\in\Lambda_{n-1}^{s-1}$ and
$\nu_j\in\Lambda_{n-1}^{s+1}$, and hence 
$$\Hom(\Delta_{n-1}(\tau),\res\Delta_n(\mu))\cong
\Hom(\Delta_{n-1}(\tau),\left(\oplus_i\Delta_{n-1}(\nu_i)\right)\oplus
\left(\oplus_j\Delta(\nu_j)\right)).$$
By semisimplicity, this Hom-space is zero unless $s+1=n-1$, i.e. unless
$s=n-2$. Thus we have reduced to considering the case $r=n$ and
$s=n-2$ as required. 
\end{proof}

Note that the test for semisimplicity in the second part of this
Theorem is typically a tractable algebra-specific calculation. This is
because for any $A_n$ satisfying (A2) (with $\lambda$ and $\mu$
as above) both $\Delta(\lambda)$
and $\Delta(\mu)$ have few composition factors (indeed the former is a
simple module). Thus the determination of homomorphisms between them
will in many cases be a tractable algebra-specific
calculation. 

It will be convenient to note the following property of algebras
satisfying (A1). Let $m<n$ with $m-n=2i$ for some $i\in\NN$. Then by
the remarks after (A1) we have that
%As $A_{n-2}\cong e_nA_ne_n$ we may identify $e_{n-2}$ with an
%idempotent (also denoted $e_{n-2}$) in $A_n$. Repeating this procedure
%we can consider for $m<n$ with $m-n$ even the (necessarily commuting)
%idempotents $e_n,e_{n-2},\ldots,e_{m+2}$ in $A_n$. Denoting the
%product of these elements by $e(n,m+2)$ we have that $e(n,m+2)$ is an
%idempotent with 
$A_m\cong e_{n,i}A_ne_{n,i}$. There is a
corresponding globalisation functor, which we denote $G_m^n$, given by
$G_m^n(N)=A_ne_{n,i}\otimes_{e_{n,i}A_ne_{n,i}}(N)$ for all
$A_m$-modules $N$. It is now an elementary exercise to verify that
\begin{equation}\label{shortG}
G_m^n(N)\cong G_{n-2}G_{n-4}\ldots G_{m}(N)
\end{equation}
for all $A_m$-modules $N$.

The value of this axiom scheme hangs on there being a large number of
concrete algebras to which it applies. We will illustrate the utility
of the theory by applying it to the contour algebra in Section
\ref{rxsec}. First though we briefly sketch some other examples of its
usefulness from the literature.

\begin{ex} (i) The Temperley-Lieb algebra. See
\cite{marbook} and \cite{cgm} for details. In this case the indexing
set is $\Lambda_n=\{n, n-2, n-4\ldots, 0\ \mbox{\rm or } 1\}$ and
$\Lambda^n=\{n\}$.  
We have a short exact sequence
$$0\too \Delta_{n-1}(i-1)\too\res \Delta_n(i)\too\Delta_{n-1}(i+1)\too 0$$
for $0\leq i<n$, and $\res\Delta_n(n)\cong \Delta_{n-1}(n- 1)$, and similar
sequences for $\ind\Delta_n(i)$.

(ii) The blob algebra $b_n$ was introduced in \cite{msblob}, and an
analysis of the form described above first carried out (in
characteristic zero) in \cite{mwblob}. These results were later
generalised to positive characteristic in \cite{cgm}.  In particular
(A1) is proved in \cite[Proposition 3]{msblob}, (A2) in
\cite[(3.2)]{mwblob}, (A3) is obvious, (A4) in \cite[Proposition
2]{msblob}, (A5) and (A6) in \cite[(3.4) Proposition and (8.2)
Theorem]{mwblob} (see also \cite[Proposition 3]{mryom}).  

In this
case the indexing set $\Lambda_n=\{n,n-2, n-4, \ldots,2-n,-n\}$ with
$\Lambda^n=\{\pm n\}$. 
We have a short exact sequence
$$0\too \Delta_{n-1}(i\mp1)\too\res \Delta_n(i)\too\Delta_{n-1}(i\pm 1)\too 0$$
for $0\leq i<n$ respectively $-n<i< 0$,  and
$\res\Delta_n(\pm n)\cong \Delta_{n-1}(\pm n\mp 1)$. There are similar
sequences for $\ind\Delta_n(i)$.

(iii) The partition algebra was introduced in \cite{mardef}. In this
case the application of the theory in this section is a little more
involved, as the tower of algebras interleaves partition algebras with
auxilliary intermediates. Details can be found in \cite{marhalf}.

(iv) Certain planar algebras --- for example planar algebras on
$1$-boxes (see \cite[Section 2.2]{jonesplan}). Planar algebras were
introduced by Jones in \cite{jonesplan} formalising and generalising
the treatment of the Temperley-Lieb algebra suggested in \cite[Section
6.2]{marbook} (and implemented in \cite{mardef} in the non-planar
setting). The verification of the axioms in this case is left as an
exercise (but see below).
\end{ex}

%}}}
%{{{ S Rui--Xi algebras\label{rxsec}
\section{Contour algebras}\label{rxsec}
%{{{ preamble + diagrams I

In this section we define a new class of algebras, the contour
algebras $X_{n,m}^d$, over a general ring $R$. We then apply the
general theory developed in the preceding section.  As we will need to
consider several different algebras, in this section we will denote
the index set for the simple modules for an algebra $A$ by
$\Lambda(A)$.

\begin{figure}[ht]
\centerline{\epsffile{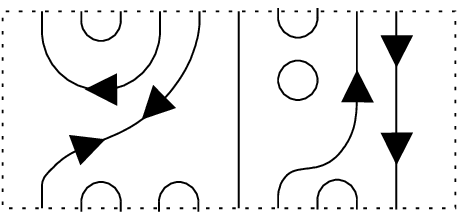}}
\caption{\label{arrow}}
\end{figure}

We will be interested in two classes of decorated Temperley-Lieb
diagrams: arrow diagrams and bead diagrams. By an {\it arrow diagram}
we mean a rectangular box containing non-intersecting line segments,
possibly with one or more arrows on each line (see Figure
\ref{arrow}).  A {\it bead diagram} is similar but with unoriented
beads instead of arrows.

It will be convenient to recall some standard terminology for ordinary
(undecorated) Temperley-Lieb diagrams which will also be needed
here. We refer to the dotted boundary of a diagram as its {\it frame}
and the interior line-segments as {\it lines}.  Lines in a diagram are
called {\it propagating lines} if they connect the northern and
southern edges of the frame, and {\it northern} (respectively {\it
southern}) {\it arcs} if they meet only the northern (respectively
southern) edge of the frame.  The endpoints of lines are called {\it
nodes}.  We identify two diagrams if they differ by an (edgewise)
frame-preserving ambient isotopy.  If the number of southern nodes in
$A$ equals the number of northern nodes in $B$ then we define the
product $AB$ to be the concatenation of the diagram $A$ above the
diagram $B$.  (In the product of two diagrams $AB$ we assume that the
southern nodes of $A$ are identified with the corresponding northern
nodes of $B$, and ignore the dotted line segment formed by their
frames across the centre of the new diagram. Then $AB$ is another
diagram.)

%}}}
%{{{ the issues figure!

%}}}
%{{{ composition -> (A1)

We say that a line in a diagram is {\it of depth $1$} (or {\it
exposed}) if the diagram can be deformed ambient isotopically such
that the line touches the eastern edge of the frame. We now define the
depth of a general line inductively by saying that a line is of {\it
depth $d$} if it is not of depth less than $d$ but can be deformed
ambient isotopically to touch a line of depth $d-1$. We say that a
diagram is {\it decorated to depth $d$} if all decorated lines in the
diagram are of depth at most $d$. For example, the diagram illustrated
in Figure \ref{arrow} is decorated to depth 5, and indeed to depth
$d$ for any $d>5$.

An arrow assigns an orientation to a line.  We say that two arrows on
the same line are {\it opposing} if they assign opposite orientations
to the line.  An arrow on a northern or southern arc is called {\it
easterly} (respectively {\it westerly}) if it point towards the
eastern (respectively western) end of the line. Similarly arrows on
propagating lines are either {\it northerly} or {\it southerly}.

Let $\bar{D}^l_n$ be the set of bead diagrams with $l$ northern and
$n$ southern nodes, and $\bar{D}_n=\bar{D}^n_n$.  The corresponding
subsets of diagrams decorated to depth $d$ will be denoted
$\bar{D}_n^l[d]$ and $\bar{D}_n[d]$ respectively.  Note that in the
composition of any two diagrams we may expose new line segments but
cannot produce new unexposed lines. Clearly similar remarks hold for
lines of depth at most $d$, and hence we have

\begin{lem}\label{contour}
The diagram product gives a map from  $\bar{D}_n^l[d]\times \bar{D}_m^n[d]$ 
to $\bar{D}_m^l[d]$.
\end{lem}

Another way to think of this is that the lines in a diagram are
contours (or isobars) and that under composition non-closed lines can
be combined to become closed contours. Fixing the eastern edge at
sea-level, the maximum physical height a contour can realise on
closure is its diagram depth. Thus depth cannot be increased by
composition.

Fix $m$, and choose elements $\delta_0,\ldots,\delta_{m-1}$ in $R$. By
Lemma \ref{contour} we may define the {\it contour algebra}
$\bar{X}_{n,m}^d=\bar{X}_{n,m}^d(\delta_0\ldots,\delta_{m-1})$ to be
the algebra obtained from $R\bar{D}_n[d]$ under concatenation with the
following additional relations:

(i) A diagram with $m$ beads on the same line is identified
with the same diagram with the beads omitted.

(ii) A diagram with an excess (modulo $m$) of $k$ beads on a given
closed loop is identified with $\delta_k$ times the same diagram with
the closed loop omitted.

\noindent It is evident that $\bar{X}_{n,m}^d$ is associative, unital,
and free as an $R$-module.

%
%By Lemma \ref{contour} we may now define the {\it contour algebra}
%$\bar{X}_{n,m}^d$ ($=\bar{X}_{n,m}^d(\delta_0,\ldots,\delta_{m-1})$)
%to be the subalgebra of $\bar{X}_{n,m}^{\infty}$ generated by those
%diagrams in $\bar{D}_n[d]$. In what follows all references to algebras
%$\bar{X}_{n,m}^d$ will include the case $d=\infty$ described above.

We denote by $\bar{X}_{n,m}^{\infty}$ the case where we allow decorated
lines of arbitrary depth. Clearly we have that $\bar{X}_{n,m}^{\infty}\cong
\bar{X}_{n,m}^n$, and for general $d$ that $\bar{X}_{n,m}^d\subseteq
\bar{X}_{n,m}^{d+1}$.

There is another presentation of these algebras in terms of arrow
diagrams. Let $D_l^n$ be the set of arrow diagrams with $l$ northern
and $n$ southern nodes, and define sets $D_n$, $D^l_n[d]$, and
$D_n[d]$ as in the corresponding bead cases. Now we define the algebra
$X_{n,m}^d$ ($=X_{n,m}^d(\delta_0,\ldots,\delta_{m-1})$) to be the
algebra obtained from $RD_n[d]$ under concatenation with the following
additional relations:

(i) A diagram with two opposing arrows on the same line is identified
with the same diagram with the two arrows omitted.

(ii) A diagram with $m$ non-opposing arrows on the same line is identified
with the same diagram with the arrows omitted.

(iii) A diagram with an excess (modulo $m$) of $k$ anti-clockwise
arrows over clockwise arrows on a given closed loop is identified with
$\delta_k$ times the same diagram with the closed loop omitted.

\begin{figure}[ht]
\centerline{\epsffile{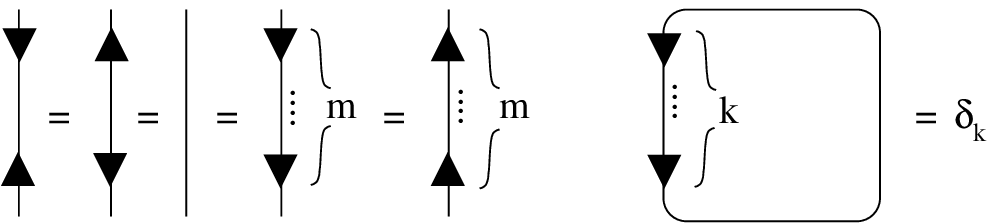}}
\caption{\label{rel1}}
\end{figure}

These three sets of relations are illustrated schematically in Figure
\ref{rel1}.

\begin{figure}[ht]
\centerline{\epsffile{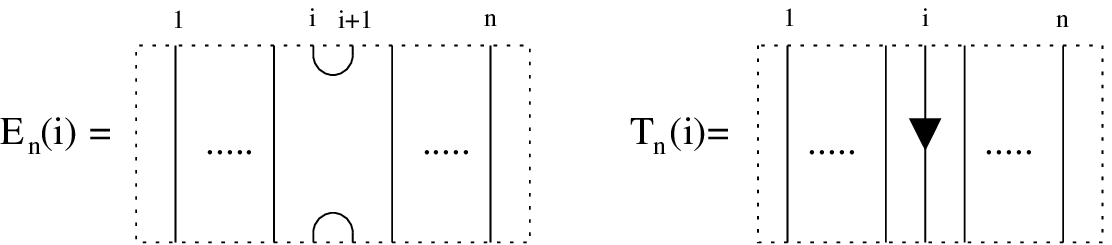}}
\caption{\label{gens}}
\end{figure}

It will be convenient to have names for certain diagrams. It is clear
that the algebra $X_{n,m}^d$ is generated by the elements $E_n(i)$ (for
$1\leq i\leq n-1$) and $T_n(i)$ (for $\max(1,n+1-d)\leq i\leq n$)
illustrated in Figure \ref{gens}. Note that
$E_n(i)^2=\delta_0E_n(i)$. The analogue of $T_n(i)$ with a bead
instead of an arrow will be denoted $\bar{T}_n(i)$.
 
\begin{prop}\label{samealg}
The map 
$$E_n(i)\longmapsto E_n(i)\quad \mbox{and}\quad T_n(i)\longmapsto
\left\{\begin{array}{ll}\bar{T}_n(i)\ &\mbox{if}\ i \ \mbox{odd}\\
(\bar{T}_n(i))^{m-1}\ &\mbox{if}\ i \ \mbox{even}\end{array}\right.$$
induces an algebra isomorphism from $X_{n,m}^d$ to $\bar{X}_{n,m}^d$.
\end{prop}
\begin{proof}
 This is an easy exercise. For example, any pair of opposing
(respectively non-opposing) arrows must arise (from these generators)
from a pair of elements $T_n(i)$ and $T_n(j)$ with $i-j$ odd
(respectively even). 
\end{proof}

Because of Proposition \ref{samealg} we will henceforth also refer to
$X_{n,m}^d$ as the contour algebra.

\begin{rem} The algebra $X_{n,m}^0$ coincides with the 
Temperley-Lieb algebra (for any $d$), while $X_{n,2}^1$ is isomorphic
to the blob algebra and $X_{n,m}^1$ to the coloured blob algebra
introduced in \cite{mwl}.  By comparing the arrow definition with that
in \cite[Definition 3.3]{rx} it is easy to show that
$X_{n,m}^{\infty}$ is isomorphic to the cyclotomic Temperley-Lieb
algebra $\tilde{T\!L}_{n,m}$ introduced by Rui and Xi (which are
planar algebras on $1$-boxes). The algebras $X_{n,m}^d$ with $1<d<n$
are new.
\end{rem}

Henceforth we  take $R=k$,
an algebraically closed field. We will show that, with some conditions
on the characteristic of $k$ and the parameters $\delta_i$, the
algebras $X_{n,m}^d$ satisfy (A1--A6).

\begin{figure}[ht]
\centerline{\epsffile{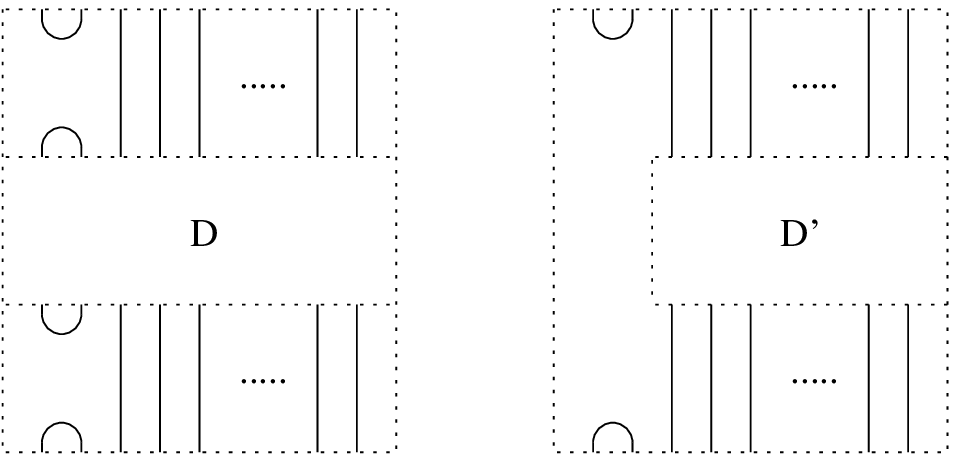}}
\caption{\label{A1fig}}
\end{figure}

\begin{prop}\label{first}
For $\delta_0 \neq 0$ we have 
$$E_n(1) X_{n,m}^d E_n(1)  \cong  X_{n-2,m}^d.$$
\end{prop}
\begin{proof}
 Any diagram $E_n(1)DE_n(1)$ in $E_n(1) X_{n,m}^d E_n(1)$ is of the
form shown on the lefthand side of Figure \ref{A1fig}, and can be put
into the form on the righthand side of the Figure for some diagram
$D'$ in $X_{n-2,m}^d$. As $\delta_0\neq 0$, the set of diagrams of the
form shown on the righthand side defines an algebra isomorphic to
$X_{n-2,m}^d$, via the map which sends $E_n(1)DE_n(1)$ to
$\delta_0D'$.
\end{proof}

This verifies (A1) when $\delta_0\neq 0$. An analogous result can be
obtained under the weaker assumption that there exists some $j$ with
$\delta_j\neq 0$. For this we argue as above, but replace every
occurrence of $E_n(1)$ with the same diagram decorated with $j$
westerly arrows on the southern arc. Henceforth we assume that there
exists some $\delta_j\neq 0$, fix $m$, and denote $X_{n,m}^d$ by $A_n$.
In proofs we will suppose that $\delta_0\neq 0$ and denote
$\delta_0^{-1}E_n(1)$ by $e_n$. The modifications for the general case
are exactly as for Proposition \ref{first} above.

%}}}
We define the {\it propagating number} of a diagram $D$ to be
the number of propagating lines in $D$. Let $D_n[d;i]$ denote the
subset of $D_n[d]$ consisting of diagrams with propagating number $i$.
Note that there is a unique undecorated diagram with no closed loops
in $D_n[d;n]$, which is the identity element in $A_n$. All other
diagrams in $D_n[d;n]$ have the same underlying undecorated diagram,
but with additional arrows and/or closed loops.  The set $D_n[d;i]$ is
not linearly independent, so we define $D^+_n[d;i]$ to be the subset
of diagrams in $D_n[d;i]$ with no closed loops, no more than $m-1$
arrows on any single line, and all arrows either westerly or
southerly. We set $D^+_n[d]$ to be the union of the $D_n^+[d;i]$. It
is easy to see that such diagrams are linearly independent, and
further that (after applying the defining relations) the composition
of diagrams restricts to a map from $D^+[d]\times D^+[d]$ to $R\times
D^+[d]$.

Let $kC_m$ be the group algebra over $k$ of the cyclic group of order
$m$.  As $T_n(i)^m=1$, the element $T_n(i)$ generates a copy of
$kC_m$.

%\noindent {\bf Remark.} 
\begin{rem}
 It is a triviality to construct an enumerated basis of
$X_{n,m}^d$ which coincides with the finite set $D_n^+[d]$, using the
technique of \cite[Proposition 2]{msblob}.  As in all the diagram
algebras mentioned in Section~\ref{towers}, this construction exhibits
bases for certain submodules of $RD_n^+[d]$ (regarded as the regular
representation).  It shows explicitly that the sum of squares of the
ranks of these submodules is the rank of $X_{n,m}^d$. These modules
coincide, in quasi-hereditary specialisations to be discussed shortly,
with the standard modules considered in Section \ref{rxrep}.  
\end{rem}

Suppose that $\delta_0\neq 0$, and consider the filtration of $A_n$ by
two-sided ideals
\begin{equation}\label{chainis}
\ldots \subset A_n E_n(1) E_n(3) A_n \subset A_n E_n(1) A_n \subset
A_n.
\end{equation}
We will denote the product $\prod_{j=1}^i E_n(2j-1)$ by
$E_{n,i}$. As $\delta_0\neq 0$ this is a preidempotent (i.e. a
non-zero scalar multiple of an idempotent), and we
define $e_{n,i}$ to be the corresponding idempotent
$\delta_0^{(-i)}E_{n,i}$. The corresponding constructions for
$\delta_j\neq 0$ are obvious.
 
\begin{prop}\label{secbas}
The $i^{th}$ section 
$$ A_ne_{n,i}A_n /A_ne_{n,i+1}A_n$$ in this filtration has basis
$D_n^+[d;n-2i]$.
\end{prop}
\begin{proof}
 This is straightforward --- confer \cite[Corollary 1.1]{msblob}.
\end{proof}

In particular we have 

\begin{cor}\label{cycis}
$$ A_n / A_n e_n A_n  \cong  (kC_m)^{\min(n,d)}   . $$
\end{cor}

A parameterisation of the simple modules of $A_n$ now follows
immediately from (\ref{simplesare}):

\begin{cor}\label{classsimple}
Suppose that there exists some $j$ with $\delta_j\neq 0$. Then for all
$n\geq 0$ we have
$$
\Lambda(X_{n,m}^d) 
\; = \; \Lambda(X_{n-2,m}^d)\sqcup(\Lambda(kC_m))^{\min(n,d)} 
\; = \;   
{\textstyle \coprod}_{i=n,n-2,..,1/0}\, (\Lambda(kC_m))^{\min(i,d)} . 
$$ 
\end{cor}

The representation theory of $(kC_m)^n $ is well understood.
For example, if $k$ is a splitting field of $x^m -1$ of characteristic
$p$ such that $p=0$ or $p$ does not divide $m$, then the set 
$ \{ 1,2,..,m \}$ may be taken as an index set $\Lambda(kC_m) $ for the
simples of $kC_m$ over $k$, and $\Lambda((kC_m)^n)=(\Lambda(kC_m))^n$. 
In the special case $d=\infty$ this provides a very short proof of
\cite[Corollary 5.4]{rx}.

Note that the restriction rules for $(kC_m)^r$ to $(kC_m)^{r-1}$ are
elementary. This will facilitate verification of (A5) shortly.

\begin{figure}[ht]
\centerline{\epsffile{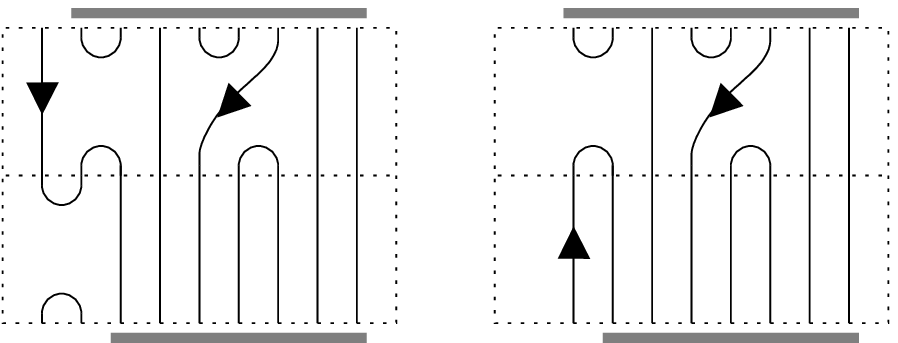}}
\caption{\label{A4fig}}
\end{figure}

Before going on to consider quasi-heredity, we quickly note that (A3)
and (A4) are both easily verified. For (A3) we can identify $A_n$
as a subalgebra of $A_{n+1}$ via the map which adds an undecorated
propagating line to the lefthand side of each diagram.  For (A4), note
that the left action of $A_{n-1}$ is by concatenation from above on
the rightmost $n-1$ strings, while the right action of $A_{n-2}$ is by
concatenation from below on the rightmost $n-2$ strings. We define a
map from a diagram in $A_ne_n$ to a diagram in $A_{n-1}$ by first
deforming the original diagram ambient isotopically to move the
leftmost northern node anticlockwise around the frame to become the
leftmost southern node, and then removing the southern arc adjacent to
this new node. An example of this is given in Figure \ref{A4fig},
where the effect of the map on the lefthand diagram is illustrated on
the right. (The shaded areas indicate the nodes acted on by the
actions from above and below.) It is easy to verify that this map
gives the desired left $A_{n-1}$-, right $A_{n-2}$- bimodule
isomorphism.

We next verify (A2).

\begin{prop}\label{qhprop}
Suppose that there exists some $j$ with $\delta_j\neq 0$, and that
either $p=0$ or $p$ does not divide $m$. Then for all $n\geq 0$ the algebra
$A_n$ is quasi-hereditary, with heredity chain of the form given in
(\ref{chainis}).
\end{prop}
\begin{proof} We consider the case $j=0$, when the heredity chain will be
precisely the chain in  (\ref{chainis}). For arbitrary $j$ we must
replace each $E_n(i)$ with the appropriately decorated analogue
introduced after Proposition \ref{first}.

We wish to show that the filtration in (\ref{chainis}) is a heredity
chain for $A_n$; i.e. that each of the quotients
$(A_ne_{n,i}A_n)/(A_ne_{n,i+1}A_n)$ is a heredity ideal of
$A_{n,i}=(A_n)/(A_ne_{n,i+1}A_n)$. For this it is enough to show that
the conditions (A2)(i) and (ii) both hold.

%The proof of (i) is entirely analogous to that of Corollary
%\ref{cycis}. By Proposition \ref{secbas} we know that
%$A_{n,i}e_{n,i}A_{n,i}$ has basis $D^+_n[d;n-2i]$. However, the subset
%of these preserved (up to scalars) in $e_{n,i}A_{n,i}e_{n,i}$ has
%all but the last $n-2i$ northern (and southern) nodes on
%non-propagating lines, and hence all remaining lines must
%propagate. Thus the quotient algebra $e_{n,i}A_{n,i}e_{n,i}$ is
%isomorphic to $(kC_m)^{\min(n-2i,d)}$, which is semisimple by our assumptions
%on $k$.
Condition (i) follows immediately from Corollary \ref{cycis} and our
assumptions on $p$.  For (ii), we begin by noting that
$A_{n,i}e_{n,i}$ has a basis represented by those diagrams with $i$
non-nested southern arcs on the $2i$ westernmost vertices, and $n-2i$
propagating lines (possibly with decorations). We have a similar basis
for $e_{n,i}A_{n,i}$ with northern instead of southern arcs. Thus the
product $D$ of such a diagram in $A_{n,i}e_{n,i}$ with such a diagram
in $e_{n,i}A_{n,i}$ must have precisely $n-2i$ propagating lines, and
it is clear that any pair of diagrams giving rise to $D$ must be
equivalent in
$A_{n,i}e_{n,i}\otimes_{e_{n,i}A_{n,i}e_{n,i}}e_{n,i}A_{n,i}$. (To see
this note that such pairs of diagrams can only differ in the
distribution of decorations between them, which can be adjusted via an
element of $e_{n,i}A_{n,i}e_{n,i}$.)

Thus we have verified (A2)(i) and (ii), and hence $A_n$ is
quasi-hereditary.
\end{proof}

%}}}
%{{{ S Representations of the Rui-Xi algebra
\section{Representations of contour algebras}\label{rxrep}

Henceforth we will assume that $A_n$ satisfies the conditions of
Proposition \ref{qhprop}. Then by the general theory in Section
\ref{towers}, every standard module $\Delta_n(\lambda)$ of $A_n$ is
the image under $G_{n-2}G_{n-4}\ldots G_{n-2i}$ of some standard
module for $(kC_m)^j$ lifted to $A_j$, for some $i,j\geq 0$ with
$2i+j=n$. (We adopt the convention that $(kC_m)^0=k$, with simple
module labelled by $\emptyset$.) We call $j$ the
propagating number of $\lambda$. Thus we need to fix our convention
for lifting modules from $(kC_m)^n$ to  $A_n$.

\begin{figure}[ht]
\centerline{\epsffile{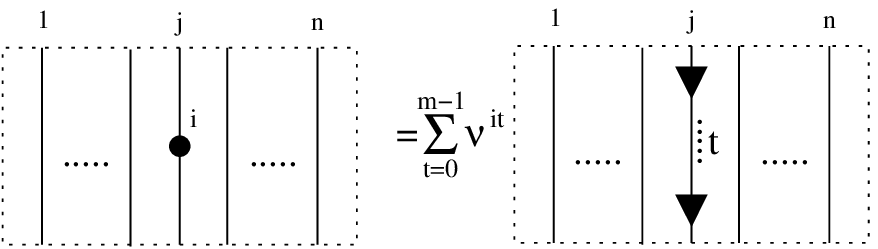}}
\caption{\label{idemis}}
\end{figure}

We fix  $\nu$, a primitive $m$th root of unity, and define
the element $\epsilon_n(i,j)=\sum_{t=0}^{m-1}\nu^{it}T_n(j)^t$ in $A_n$
(where $T_n(j)^0=1_{A_n}$). Note that this element is a preidempotent: we
have $(m^{-1}\epsilon_n(i,j))^2=m^{-1}\epsilon_n(i,j)$. Graphically we
represent $\epsilon_n(i,j)$ as shown in Figure \ref{idemis} and refer
to its decoration as $\bullet(i)$.

Now the simple module labelled by $(i_1,\ldots, i_n)$ for $(kC_m)^n$
can be realised as an $A_n$-module (via Corollary \ref{cycis}) as the
module $A_n\epsilon_n(i_1,1)\ldots \epsilon_n(i_n,n)$, with the
convention that we identify any diagram with fewer than $n$
propagating lines with zero. There is an obvious extension of the
graphical notation for $\epsilon_n(i,j)$, where we represent
$\epsilon_n(i_1,1)\ldots \epsilon_n(i_n,n)$ by the corresponding
product of the diagrams for each $\epsilon_n(i,j)$.

\newcommand{\lm}{l}
By the general theory in Section \ref{towers} we have for $n>\lm$ with
$n-\lm$ even that 
$$
\Delta_n(i_1,\ldots,i_{\lm})
\cong G_{\lm}^n\Delta_{\lm}(i_1,\ldots,i_{\lm}) \cong A_ne_{n,t}
\otimes_{e_{n,t}A_ne_{n,t}}A_{\lm} \epsilon_{\lm}(i_1,1)\ldots
\epsilon_{\lm}(i_{\lm},\lm)
$$ where $t=\frac{n-\lm}{2}$.  Let $D^n_{\lm}(i_1,\ldots,i_{\lm})$
denote the set of diagrams with $n$ northern and $\lm$ southern nodes,
$\lm$ propagating lines and no closed loops, such that the $j$th
propagating line is decorated with $\bullet(i_j)$.  Let
$\Delta'_n(i_1,\ldots,i_{\lm})$ denote the $A_n$-module with basis
$D^n_{\lm}(i_1,\ldots,i_{\lm})$, where the action of $A_n$ is by
concatenation from above, such that any product of diagrams with fewer
than $\lm$ propagating lines is set to zero.  It will be evident that
a fixed distribution of southern arcs could be added to every diagram
without changing the action, and hence we have

\begin{prop}\label{stanbasis}
The modules $\Delta_n(i_1,\ldots,i_{\lm})$ and
$\Delta'_n(i_1,\ldots,i_{\lm})$  can be identified. 
\end{prop}

We now consider (A5) and (A6). We first note that there is an 
$A_{n-1}$-submodule of $\Delta_n(i_1,\ldots,i_{\lm})$ (as a diagram module)
spanned by those diagrams with a propagating line from the most westerly
northern node is isomorphic to $\Delta_{n-1}(\mu)$ where
$\mu=(i_2,\ldots,i_{\lm})\in \Lambda_{n-1}^{\lm-1}$. (This is clear, as
$A_{n-1}$ acts on all but the most westerly northern node.)

\begin{figure}[ht]
\centerline{\epsffile{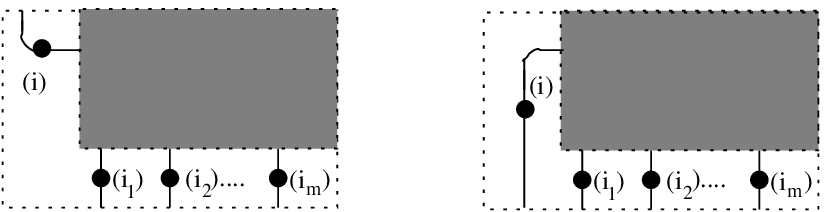}}
\caption{\label{rest}}
\end{figure}

All remaining diagrams in $\Delta_n(i_1,\ldots,i_{\lm})$ have a
northern arc starting at the most westerly northern node. We consider
a new basis for this set formed by taking linear combinations of
diagrams such that this northern arc is decorated with a $\bullet(i)$
for some $i$, as illustrated in the left-hand diagram in Figure
\ref{rest} (where the shaded region denotes some collection of lines
whose precise configuration does not concern us). If we take the
subset of such diagrams with fixed decoration $\bullet(i)$ then,
modulo the submodule $\Delta_{n-1}(\mu)$ described above, there is an
$A_{n-1}$-module isomorphism with $\Delta_{n-1}(\nu)$ (where
$\nu=(i,i_1,\ldots,i_{\lm})\in\Lambda_{n-1}^{\lm+1}$) given by the map
which deforms the diagram ambient isotopically as shown in Figure
\ref{rest}.

This completes the verification of (A5); it is also clear from the
above that (A6$'$) holds. Thus we may apply all the general theory
from Section \ref{towers} to these algebras.

To apply Theorem \ref{homcatch} it only remains to calculate
$\Hom(\Delta_n(\lambda),\Delta_n(\mu))$ for all
$\lambda\in\Lambda_n^n$ and $\mu\in\Lambda_n^{n-2}$. If there exists a
$\mu\in\Lambda_n^{n-2}$ with $\Delta_n(\mu)$ non-simple, then at least
one such Hom-space will be non-zero. Thus to prove that our algebras
are semisimple it is enough, for example, to show that the Gram matrix
for $\Delta_n(\mu)$ is non-degenerate for all $\mu\in\Lambda_n^{n-2}$.
%We do this, following some explicit examples, in the next section
%(Proposition~\ref{grammy}) for generic values of the parameters.
%\footnote{These are precisely the Gram matrices
%calculated in \cite[Proposition 8.1]{rx}. The answer given there is an
%explicit but complicated non-zero polynomial in the defining
%parameters.} 
%Subject to this we can deduce
%
%\begin{thm}\label{genssim} 
%The algebras $A_n$ are generically semisimple.
%\end{thm}\qed

%}}}
%{{{ Gram examples
\section{Gram matrix results}\label{gramres}

We now consider the Gram matrix $G_n(\lambda)$ of inner products with
respect to the diagram basis of $\Delta_n(\lambda)$ (confer
\cite{msblob}).  \newcommand{\plam}[1]{i_1,\ldots,i_{#1}} Let
$D^n_{l,p}(\plam{p})$ be the mild generalisation of $D^n_p(\plam{p})$
consisting of diagrams with $n$ northern and $l$ southern nodes, 
$p$ propagating lines, and propagating line decorations as for
$D^n_p(\plam{p})$.  Then $\hat{D}^p_n(\plam{p})=D^p_{n,p}(\plam{p})$ is
the upside down version of $D^n_p(\plam{p})$.
%and hence a basis for
%the dual of $\Delta_n(\lambda)$ (where we set $\lambda=(\plam{p})$).
Let $\epsilon(\lambda)$ denote the unique element of $D^p_p(\lambda)$.
Consider the map
\begin{eqnarray*}
\hat{D}^p_n(\lambda) \times D^n_p (\lambda) 
 & \rightarrow  & {\mathbb Z}[\delta_0,\ldots,\delta_{m-1}]
\\
(a,b) & \mapsto & \langle a | b \rangle
\end{eqnarray*}
where $\langle a | b \rangle$ is such that the diagram product $ab \;
= \; \langle a | b \rangle \; \epsilon(\lambda)$ if $ab$ lies
in ${\mathbb Z}[\delta_0,\ldots,\delta_{m-1}] D^p_p(\lambda)$, and is
zero otherwise.  Note that $\langle - |- \rangle$ defines an inner
product on $\Delta_n(\lambda)$.
%{{{ m=2 and first example

We will first consider the case $m=2$ and $d=\infty$ for the sake of
definiteness.  However, neither restriction is significant. In
pictures we will denote $\bullet(1)$ just by $\bullet$.  When $n=2$ we
then have
$$
\Lambda_2 = \Lambda^2 \cup \Lambda_0 
= \{ (1,1), (1,2), (2,1), (2,2) \} \cup \{ \emptyset \}
$$
(using the index set introduced above Corollary~\ref{classsimple}). 

\begin{figure}[ht]
\centerline{\epsffile{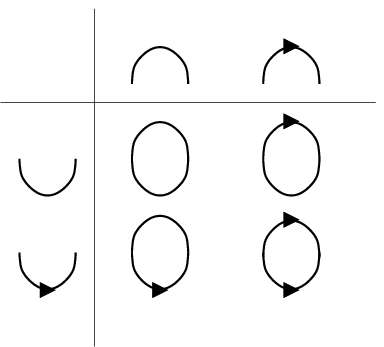}}
\caption{\label{ex1}}
\end{figure}

%%It follows from the relations (i)--(iii) on diagrams above
A simple restatement of the inner product above is that we need only
 consider the concatenation of the top halves of diagrams in the
 diagram basis of a standard module with bottom halves in the
 dual. Accordingly we may compute the Gram matrix $G_2(\lambda)$ for
 $\Delta_2(\lambda)$ with $\lambda=\emptyset$ from the diagrams in
 Figure \ref{ex1}, which give the corresponding matrix
$$\left( \begin{array}{cc} \delta_0 & \delta_1 \\ \delta_1 & \delta_0
\end{array} \right).$$
That is, $|G_2(\emptyset)|=\delta_0^2 -\delta_1^2$. 

%\[
%\epsffile{mn22.eps}
%\hspace{30pt} \vspace{-30pt}
%= 
%\vspace{30pt}
%\]

%}}}
%{{{ singular moment

Let us consider for a moment what happens in a singular
specialisation. If $\delta_1 = \delta_0$, then $\Delta_2(\emptyset)$
is not simple.  Armed with this knowledge it is straightforward to
construct a proper submodule.  Indeed it will be evident that if we
write $a$ and $b$ for the two basis elements depicted, then $T_2(i)
(a-b) = -(a-b)$ for $i=1,2$, and $E_2(1) (a-b) = 0$.  Thus $(a-b)$
generates a submodule of $\Delta_2(\emptyset)$ isomorphic to
$\Delta_2(1,1)$ in such a specialisation.  
%The resultant standard
%module homomorphism engenders a swathe of homomorphisms throughout the
%family (following the arguments in the proof of
%Theorem~\ref{homcatch}). For example there is a homomorphism
By Theorem \ref{homcatch}(i) we obtain corresponding homomorphisms
$$\Delta_n(1,1) \rightarrow \Delta_n(\emptyset)$$
for all even $n$. 

%}}}
%{{{ next example+conclusions

\begin{figure}[ht]
\centerline{\epsffile{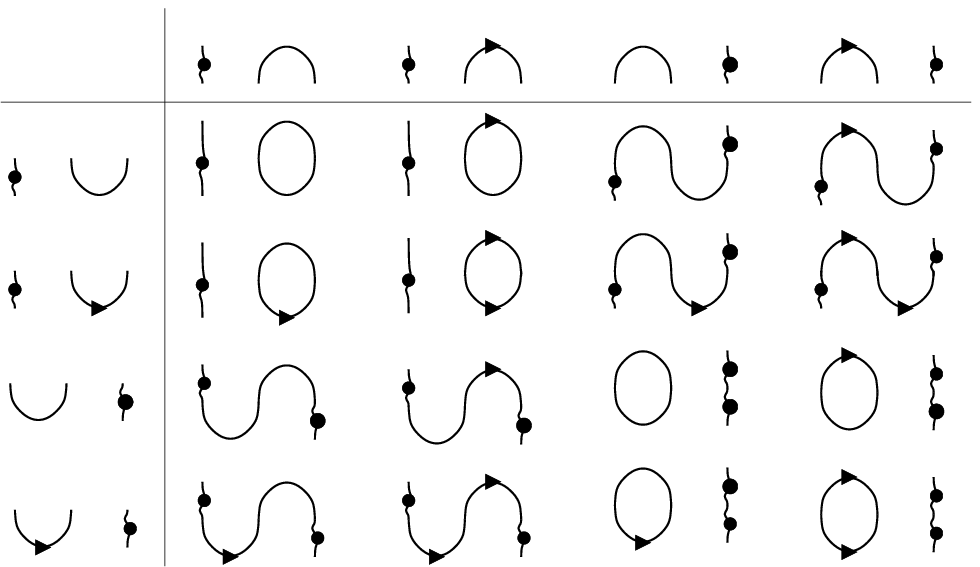}}
\caption{\label{ex2}}
\end{figure} 

Returning to generic parameters, 
for $\Delta_3(\lambda)$ with $\lambda=(1)$ or $(2)$ we have from
Figure \ref{ex2} that the Gram matrix equals
$$\left( \begin{array}{cccc} 
\delta_0 & \delta_1 & 1 & \pm 1 \\ 
\delta_1 & \delta_0 & \pm 1 & 1 \\
 1 & \pm 1 & \delta_0 & \delta_1 \\
 \pm 1 & 1 & \delta_1 & \delta_0
\end{array} \right)$$

%\epsffile{mn23.eps}
%\hspace{30pt} %\vspace{-30pt}
%= 
%\]

The determinant here is again easy to compute, but the details do not
concern us here. Instead we return to the general case.

\begin{prop}\label{grammy}
Considering $\delta_0,\delta_1,\ldots,\delta_{m-1}$ as indeterminates, 
the determinant $| G_n(\lambda) |$ is non-zero.
\end{prop}
\begin{proof}  It is clear that all Gram matrix elements take the form $\xi
\prod_i (\delta_i)^{\alpha_i}$ where $\xi$ is some $m$th root of
unity.  Consider for a moment the diagonal elements of the Gram
matrix, organised as indicated by our examples.  In these, every upper
arc meets a mirror image lower arc, and either both are undecorated,
or they have `cancelling' decorations.  Thus every arc contributes
positively to $\alpha_0$.  It follows that in each row of any Gram
matrix the value of $\alpha_0$ for the matrix element on the diagonal
strictly exceeds any other, and hence that $| G_n(\lambda) |$ is a
non-zero polynomial.
\end{proof}

\begin{cor} 
The algebras $X_{n,m}^d$ are generically semisimple with respect to
the Zariski topology for our parameter space.
\end{cor}

%}}}
%}}}
\section{Discussion}\label{discus}

Note that we have just proved generic semisimplicity of our algebras
without appeal to the full strength of the machinery developed in
Section \ref{towers}. However, Proposition \ref{grammy} does not
provide a means for determining {\it which} specialisations are
non-semisimple; indeed determining the zeros of $| G_n(\lambda) |$ for
general $\lambda$ seems a rather intractable problem. We conclude by
discussing how our result can be strengthened using the machinery
developed.

By Theorem \ref{homcatch}(ii), we have the much simpler condition

\begin{cor} \label{simplertest}
The algebra $X_{n,m}^d$ is semisimple over $k$ if and only
 if $(\delta_0,\ldots,\delta_{m-1})$ is such that 
$$\prod_{n'\leq n}\prod_{\lambda\in\Lambda_{n'}^{n'-2}}|G_{n'}(\lambda) 
|\neq 0.$$ 
\end{cor}

\begin{rem}For
$X_{n,m}^{\infty}$ the Gram matrices in Corollary \ref{simplertest}
are precisely those calculated in \cite[Proposition 8.1]{rx}. The
answer given there is a complicated but explicit polynomial in the
defining parameters. Thus, using the polynomial in \cite[Proposition
8.1]{rx}, we can determine precisely which specialisations of
$X_{n,m}^{\infty}$ are semisimple. Very similar explicit results may
be obtained for the algebras $X_{n,m}^d$; for $d=0$ these were
calculated in \cite{marbook}, and for $d=1$ in \cite{msblob}.
\end{rem}

The theory developed in Section \ref{towers} also provides a means for
studying non-semisimple specialisations, as it provides a means for
determining a large number of homomorphisms.  In the interests of
brevity we do not pursue the structure of the non-semisimple cases of
the  contour algebras further here.  Note, however, that much (in some
cases essentially all) of the structure of the other algebras
mentioned in Section~\ref{towers} has been derived in the literature
using methods which are entirely based on (ad hoc formulations of)
(A1-6). Similar efficacy may be anticipated here.

The second author and Ryom-Hansen recently made play with an
interesting tensor space representation of the blob algebra, which
they show in \cite{mryom} to be a full tilting module in
quasi-hereditary specialisations of that algebra. It is worth noting
that the bulk of the machinery they use in their proof follows from
our (A1-6).

In particular, suppose that we have a tower of algebras $A_n$
satisfying (A1-6), together with a contravariant duality $^o$ on each
$A_n$. For each $n$ let $T_n$ be an $A_n$-module such that

\medskip
\begin{list}{\quad \quad{\bf (A\arabic{axiom})}}{\usecounter{axiom}}
\setcounter{axiom}{6}
\item {\it (i) $T_0$ and $T_1$ are tilting modules.

$\phantom{\  \quad\quad}$(ii) For each $n\geq 2$ we have $F_n(T_n)\cong T_{n-2}$ and
  $T_n^o\cong T_n$.

$\phantom{\quad\quad}$(iii) The natural map $G_{n-2}F_n(T_n)\too T_n$ is injective.}
\end{list} 
\medskip

Then by the results in \cite[Proposition 5]{mryom} we have that $T_n$ is a
tilting module for each $n$.

Diagram algebras typically have a contravariant duality given by
inverting the individual diagrams. Thus the examples discussed in
Section 1 (together with the contour algebras) do satisfy the
conditions before (A7).  In many examples modules satisfying (A7)
arise by constructing analogues of \lq tensor space\rq\
representations for the corresponding families of algebras. We do not
have a candidate for a full tilting module here, but if one were
forthcoming then a similar analysis should be possible.

Note that the contour algebras can be further generalised by allowing
diagrams to have more than one line from a given node and/or dropping
the non-crossing condition. An obvious example would be a decorated
version of the partition algebra. The notion of depth is no longer
meaningful, and the proof of quasi-heredity is slightly more
complicated, but otherwise our machinery continues to apply. The most
significant complication is the replacement of the cyclic group in our
analysis by other, more complicated, group algebras.

We conclude with some remarks on our choice of axiom scheme. In (A1),
the choice of $N=2$ in the definition of $\Phi:A_{n-N}\longrightarrow
A_n$ could be varied.  However, for larger values of $N$ the analysis
of the interplay between induction/restriction and
globalisation/localisation becomes more complicated, and the case
$N=2$ seems to cover all diagram algebra examples introduced to
date. The reason for having intermediate layers is to ensure that
$\Delta$-restriction is multiplicity free --- a useful feature in
practical calculations (see \cite{verokun}).

Note that the heredity chain for any quasi-hereditary algebra gives
rise to a tower satisfying (A1) and (A2). It is the extra structure
imposed by the remaining axioms that we wish to emphasise here.  In
particular the metric structure induced on our set of weights by the
local behaviour (A5) justifies the use of the term weights, by analogy
with \cite{jantzen}.

Quasi-heredity is quite a strong property for an algebra to possess,
and there have been several alternatives proposed for the study of
wider classes of algebras. Important examples are cellular algebras
\cite{gl} (but see also \cite{kxcell}), tabular algebras
\cite{greentab}, and various types of stratified algebras
\cite{cpsstrat, dstrat}. It would be interesting to consider how axiom
(A2) might be weakened in these (or other) settings.

\begin{acknowledgements}
The authors are grateful to Volodymyr Mazorchuk for several useful
suggestions and references.
\end{acknowledgements}

\bibliographystyle{amsalpha} 
\bibliography{/home/anton/Work/Lib/books}
\end{document}